\documentclass{commat}

\usepackage[pdftex]{graphicx}

\newcommand{\R}{\mathbb{R}}
\newcommand{\N}{\mathbb{N}}
\newcommand{\C}{\mathbb{C}}
\newcommand{\B}{\mathbb{B}}
\newcommand{\Z}{\mathbb{Z}}
\newcommand{\ZG}{\mathbf{Z}}
\newcommand{\D}{\mathbb{D}}
\newcommand{\Ng}{\mathbf{N}}
\newcommand{\Bg}{\mathbf{B}}
\newcommand{\SG}{\mathfrak{S}}
\newcommand{\Rj}{\mathrm{Re}}
\newcommand{\Ij}{\mathrm{Im}}
\newcommand{\eg}{\mathbf{e}}
\newcommand{\ug}{\mathbf{u}}
\newcommand{\ig}{\mathbf{i}}
\newcommand{\jg}{\mathbf{j}}
\newcommand{\kg}{\mathbf{k}}
\newcommand{\Ar}{\mathrm{Arg}}
\newcommand{\ar}{\mathrm{arg}}

\title{%
   Bicomplex numbers as a normal complexified $f$-algebra
}

\author{%
     Hichem Gargoubi and Sayed Kossentini
}

\affiliation{
    \address{Hichem Gargoubi --
    Universit\'e de Tunis, I.P.E.I.T., Department of Mathematics, 2 Rue Jawaher Lel Nehru, Monfleury, Tunis, 1008 Tunisia.
}
    \email{%
    hichem.gargoubi@ipeit.rnu.tn
}
    \address{Sayed Kossentini --
    Universit\'e de Tunis El Manar, Facult\'e des Sciences de Tunis, Department of Mathematics, 2092, Tunis, Tunisia.
}
    \email{%
    kossentinisayed@hotmail.fr
}
}

\abstract{%
    The algebra $\B$ of bicomplex numbers is viewed as a complexification of the Archimedean $f$-algebra of hyperbolic numbers $\D$. This lattice-theoretic approach allows us to establish new properties of the so-called $\D$-norms. In particular, we show that $\D$-norms generate the same topology in $\B$.  We develop the $\D$-trigonometric form of a bicomplex number which leads us to a geometric interpretation of the $n^{th}$ roots of a bicomplex number in terms of polyhedral tori.  We use the concepts developed, in particular that of Riesz subnorm of a $\D$-norm, to study the uniform convergence of the bicomplex zeta and gamma functions.  The main result of this paper is the generalization to the bicomplex case of the Riemann functional equation and Euler's reflection formula.
}



\dedication{In memoriam Christian Duval}

\firstpage{125}

\VOLUME{30}

\DOI{https://doi.org/10.46298/cm.9312}

\begin{paper}

It has been shown in a recent paper~\cite{HS} that the unique (up to isomorphism) algebra that is both Clifford algebra and Archimedean $f$-algebra containing $\R$ is the algebra of hyperbolic numbers
\begin{align*}
	\D =	\left\lbrace z = x+\jg y:\quad x,y \in \mathbb{R},\, \jg \notin\mathbb{R}; \jg^{2} = 1 \right\rbrace.
\end{align*}
This noteworthy fact connects two classical domains of mathematics: Clifford algebras and Riesz spaces.
The aim of this paper is to study
the complexification of the algebra $\D$ in the framework of Riesz space theory.

According to Arnold~\cite{KT}, \textit{attempts to complexify and to quaternionize mathematical theories are making clear the fundamental unity of all parts of mathematics\dots Complexifica\-tion is an informal operation for which there are
no axioms; we should try to guess}.

In the same vein, the natural question seems to be
the following: is the algebra $\B = \D + \mathbf{i}\D$ of bicomplex numbers
a simple "multiplication" of the algebra $\D$ or is it rather a "good" complexification of the structure of Archimedean $f$-algebra of $\D$?
In this paper we will give elements of answer to this general question.
In particular, we prove that $\B$ is a normal complexified $f$-algebra.
This theoretic-lattice consideration allows us to establish
$\D$-extensions of well-known properties
of complex numbers to $\B$.

Recall that the algebra of bicomplex numbers
\begin{equation*}
\mathbb{B}: = \Big\{a+b\mathbf{i}+c\mathbf{j}+d \mathbf{k}:\, a,b,c,d \in \mathbb{R};\,
\mathbf{i}^{2} = \mathbf{k}^{2} = -1,\,\mathbf{j}^{2} = 1,\,\mathbf{i}\mathbf{k} = \mathbf{k}\mathbf{i} = \mathbf{j} \},
\end{equation*}
was
introduced in 1892 by Segre~\cite{SER} in his search for special algebras,
and since then there was considerable activity in the field for several years. One can cite, for instance, the paper of
Scorza Dragoni~\cite{DRAG} in 1934
on holomorphic functions of a bicomplex variable, the work of
Morin in 1935 on the algebra of bicomplex numbers~\cite{MOR}, the series of papers by Spampinato in 1935 and 1936 (see~\cite{SP1}, \cite{SP2}, \cite{SP3}) on functions of a
bicomplex variable, and the development of a generalized bicomplex variable due to Takasu~\cite{TAKA} in 1943.

Bicomplex numbers has found applications in geometry and quantum physics (see e.g.\cite{Ham}, \cite{Krav}) as a commutative four dimensional algebra that
generalizes complex numbers. In fact, the algebra of bicomplex numbers is the unique commutative complex Clifford algebra that is not a division algebra~\cite{BDS}:
$\B\cong \mathrm{Cl}_{1}(\C)$ which has complex numbers $\C\cong \mathrm{Cl}_\R(0,1)$ and hyperbolic numbers $\D\cong \mathrm{Cl}_\R(1,0)$ as Clifford subalgebras~\cite{SO}.

The research on bicomplex and hyperbolic numbers has been revived some
decades ago by Yaglom~\cite{Yag}, Riley~\cite{RIL}, and Price~\cite{Price}.
 In recent years, several properties of complex analysis have been generalized for bicomplex numbers. The bicomplex Riemann zeta function is introduced by Rochon in~\cite{ZET} and bicomplex quantum mechanics is investigated in~\cite{DRQ1}, \cite{DRQ2}. Kumar et al. studied bicomplex $C^\ast$-algebras and topological bicomplex modules
in~\cite{C*alg} and in~\cite{TOPB} respectively. Hahn-Banach theorems for bicomplex modules have been proved by Luna-Elizarrar\'as et al. in~\cite{HAN}. Alpay et al.~\cite{SHAP} developed the functional analysis with bicomplex scalars. Further applications and properties of bicomplex numbers can be found in~\cite{COLOM1}, \cite{COLOM2},~\cite{GERV},~\cite{LUNA1}.

It is well known that quaternions
introduced in 1843 by Hamilton~\cite{Hamil} are the only possible four-dimensional real algebra without zero divisors. Quaternions form a field but are not commutative.
From a purely algebraic point of view, the loss of commutativity is not such a big problem, but it prevents from developing a viable~\footnote{Regardless of the existence of several successful theories on holomorphicity in the quaternionic sitting such as the theory of "regular functions" initiated by Fueter in 1936.} theory of holomorphic functions of a quaternion variable.
In return, and despite the existence of zero divisors, many authors agree (see e.g.~\cite{SHAP},~\cite{LUNA1}) that bicomplex numbers can represent a reasonable alternative to quaternions to build a theory of functions of several complex variables.

In the present paper we consider a new direction. It consists in looking at the algebra $\D$, somehow, as an intermediate object between $\R$ and $\C$ (we show in an upcoming article that it is "closer" to $\R$ than to $\C$).
We believe that the fundamental structural difference between $\R$ and $\C$ is that of Archimedean $f$-algebra ($\C$ can not be endowed with such a structure~\cite{HS}). Therefore, to extend complex analysis in general
to higher dimensions with an underlying order structure, the key idea is to extend, in a manner to define, the structure of Archimedean $f$-algebra. The complexification of $\D$ can be seen as the first step in this direction.

It is clear that the difficulty in using bicomplex numbers is that $\B$ is not a division algebra. As a consequence, one crucial difference with the complex case is that in $\B$ there is no multiplicative norm. Therefore, we introduced (Proposition~\ref{jnorme}) what we called a Riesz subnorm of a $\D$-norm, which is submultiplicative, to make the algebra $\B$ into a real Banach algebra.

Our goal in this paper is twofold. The first, is to give several new concepts of bicomplex analysis and geometry based on the structure of normal complexified $f$-algebra of $\B$. After a brief reminder (Section~2) of basic properties of bicomplex and hyperbolic numbers necessary for this article, we introduce in section 3 the notion of $\D$-trigonometric form of bicomplex numbers and some of their properties. For example, we prove that the bicomplex $n^{th}$ roots of unity can be represented by the vertices of a regular polyhedral torus. In section 4 we introduce the notion of $\D$-norm on bicomplex numbers. In particular, we define the notion of Riesz subnorm of a $\D$-norm.

The second goal is to use the obtained lattice-theoretical results to go further in the development of the theory of bicomplex
zeta function introduced by Rochon in~\cite{ZET}. We establish uniform convergence of the bicomplex Riemann zeta function and define the bicomplex gamma function as an absolute convergent integral. Furthermore, a bicomplex Mellin integral and functional equations are obtained.

 The main result of this paper is the following theorem, known as bicomplex Riemann functional equation and Euler's reflection formula.

\begin{theorem}\label{main}
The following statements are satisfied.
\begin{itemize}
\item [$(\rm{i}) $] $\displaystyle\Gamma(1+\omega) = \omega \Gamma(\omega)\quad \mbox{for~} \omega \in \Omega_- $;
\item [$(\rm{ii}) $] $\displaystyle \Gamma(1-\omega) \Gamma(\omega) = \frac{\pi}{\sin \pi \omega} \quad \mbox{for~} \omega \in \Omega$;
\item [$(\rm{iii}) $] $\displaystyle\zeta(\omega) = 2 {(2 \pi)}^{\omega-1} \sin (\frac{\pi}{2} \omega) \Gamma(1-\omega) \zeta(1-\omega) \quad\mbox{for~} \omega \in 1+\B_\ast.$
\end{itemize}
\end{theorem}

The proof of Theorem~\ref{main} and details of the notations will be given in section~5.

\section{Preliminaries}

In this section we present some basic properties of hyperbolic numbers and bicomplex numbers. For more details see~\cite{SHAP},~\cite{HS}, \cite{LUNA2},~\cite{DRA},~\cite{SR} and~\cite{SO}. For the used lattice concepts we refer the reader to the standard books~\cite{Lux} and~\cite{ZAN}.

\subsection{Riesz spaces and $f$-algebras}
An ordered real vector space $L$ is said to be Riesz space or (vector lattice) if the supremum $u\vee v$; equivalently, the infimum $u\wedge v$ of two elements $u$ and $v$ exist in $L$. In this case the absolute value of $u \in L$ is defined by $|u| = u\vee(-u)$.

A Riesz space $L$ is said to be Archimedean if $\inf\{u n^{-1}:n = 1, \cdots\} = 0$ for all $u\in L^+$, where $L^+$ is the set of all positive elements of $L$ called the positive cone of $L$. A real algebra $A$ (associative algebra with usual algebraic operations) is said to be an $f$-algebra if $A$ is a vector lattice in which the positive cone $A^+$ satisfies the properties: $ a, b \in A^+$ then $a b \in A^+$; $a\wedge b = 0$ implies $ac \wedge b = a\wedge cb = 0$ for all $c \in A^+$. In any $f$-algebra the squares are positive and the absolute value is multiplicative.

A typical example of $f$-algebras is the linear space of real valued continuous functions on a topological space. Moreover, Archimedean $f$-algebras are known to be commutative (see e.g.\cite{HP}) and are even automatically associative~\cite{BE}. Of course, the fundamental example of Archimedean $f$-algebras is the field $\R$ of real numbers.

\subsection{Hyperbolic numbers}
In this section we present the results of~\cite{HS} used throughout this paper.

The ring of hyperbolic numbers
\begin{equation*}
\D: = \Big\{z = x+\jg y:\, x,y \in \R,\, \jg \notin \R;\, \jg^2 = 1\Big\},
\end{equation*}
defined in the introduction has zero divisors which are the multiples $x(1\pm \jg)$ with $x \in \R\setminus \{0\}$. Thus, the group $\D_\ast$ of units of $\D$ is characterized by all hyperbolic numbers $z$ such that $\|z\|_h\neq 0$ where $\|z\|_h$ is the \textit{hyperbolic square modulus} of $z = x+ \jg y$ defined by
\begin{equation*}
\|z\|_h: = z \bar{z} = x^2 -y^2,
\end{equation*}
where $\bar{z}$ is the conjugate of $z$ given by $\bar{z} = x-\jg y$.

The hyperbolic plane has an important basis defined by the two idempotent elements
\begin{equation*}
\displaystyle \eg_1: = \frac{1+\jg}{2} \mbox{~and~}\eg_2: = \frac{1-\jg}{2}\Longrightarrow \eg_1+ \eg_2 = 1,\, \eg_1 \eg_2 = 0.
\end{equation*}
It follows that, each hyperbolic number $z$ can be expressed uniquely as
\begin{equation}\label{specdecomp}
z = \pi_1(z) \eg_1+\pi_2(z) \eg_2,
\end{equation}
where $ \pi_1(x+\jg y) = x+y$ and $\pi_2(x+\jg y) = x-y$. The representation~\eqref{specdecomp}, called \textit{spectral decomposition}~\cite{SO}, allows us to reduce algebraic operations into component-wise operations. Moreover, the partial order defined by

\begin{equation*}
	z, w \in \D;\, z\leq w \mbox{~if and only if~}\pi_{k}(z)\leq \pi_{k}(w),\quad (k = 1,2),
\end{equation*}
	makes $\D$ into Archimedean $f$-algebra where the lattice operations are given by
\begin{equation}\label{eqsup}
z\vee w = \max \left \{\pi_1(z), \pi_1(w)\right\} \eg_1+\max \left\{\pi_2(z), \pi_2(w)\right\} \eg_2,
\end{equation}
\begin{equation}\label{eqinf}
z\wedge w = \min \left\{\pi_1(z), \pi_1(w)\right\} \eg_1+\min \left\{\pi_2(z), \pi_2(w)\right\} \eg_2.
\end{equation}
Note that the set of positive hyperbolic numbers $\D^+$, and therefore ordering in $\D$ was introduced first in~\cite[Section 2]{GERV} and considered in~\cite[Section 1.4]{SHAP} with the aim to generalize usual concept of real norm.

The Riesz space $\D$ is \textit{Dedekind complete}, that is, every nonempty set of $\D$ that is bounded from above (resp. from below) has a supremum (resp. a infimum).

For $z,w \in \D$ write: $z<v$ when $(w-z) \in \D^+ \setminus \{0\}$ and $z\ll w$ when $(w-z) \in \D_\ast^+$. So that, $z,w \in \R$ implies $z< w$ in $\R$ if and only if $z\ll w$ in $\D$.

Let $a, b \in \D$ be such that $a\leq b$. The set
\begin{equation*}
\left[a,b \right]_\D = \{z \in \D:\, a\leq z \leq b\}
\end{equation*}
is called hyperbolic closed interval. Similarly, one can define open interval $\left(a,b\right)_\D$ or semi-open intervals $\left(a,b \right]_\D $ and $\left[a,b\right)_\D $, replacing $\leq$ by $\ll$ in left-right and left/right, respectively.

From~\eqref{eqsup} the absolute value of an hyperbolic number $z$ is given by
\begin{equation}\label{Abs}
|z|: = z\vee(-z) = |\pi_1(z)|\eg_1+|\pi_2(z)| \eg_2.
\end{equation}
Thus, the kernel of group homomorphism $ |.|$ from $\D_\ast$ to $\D_\ast^+$ is the four Klein group $\SG = \{1,-1, \jg, -\jg \}$ called \textit{group of signs} of $\D$. For $\varepsilon \in \SG$, the set $\D^\varepsilon: = \{z\in \D: |z| = \varepsilon z\}$ is called the \textit{$\varepsilon$-cone} of $\D$. The $(1)$-cone is the positive cone $\D^+$ and the $(-1)$-cone is the negative cone $\D^-$. The absolute value function yields a norm in $\D$ given by the formula
\begin{equation}\label{riesznorm}
\|z\|_{R}: = \min \Big\{\alpha \in \R^+ :\, \alpha.1\geq |z|\Big\} = |z|\vee \overline{|z|} \quad \mbox{for all~} z\in \D,
\end{equation}
 and satisfying the following properties for all $z, w \in \D$:
\begin{itemize}
	\item [ $N1)$] $\|z\|_R\leq \|w\|_R$ whenever $|z|\leq |w|$;
	\item [ $N2) $] $\|1\|_{R} = 1$, $\|z w\|_{R}\leq \|z\|_{R} \| w\|_{R}$.
\end{itemize}
It follows from the properties above that $(\D, \|.\|_R)$ is a Banach lattice. Consequently, the exponential of any hyperbolic number $z$ can be defined by the absolute convergent series
\begin{equation*}\label{expoh}
e^z : = \displaystyle\sum_{n = 0}^{\infty} \frac{z^n}{n!} = e^{\pi_1(z)} \eg_1+ e^{\pi_2(z)} \eg_2.
\end{equation*}
According to the above spectral decomposition of $e^z$ one can easily verify that the hyperbolic exponential function $\mathrm{\exp}$ is a group isomorphism from $\D$ to $ \D_\ast^+$ that preserves conjugation and lattice operations :
\begin{equation}\label{expop}
\overline{e^{z}} = e^{\bar{z}};\quad e^{z}\wedge e^{w} = e^{z\wedge w}; \quad e^{z}\vee e^{w} = e^{z\vee w}\quad\mbox{for all~} z,w \in \D.
\end{equation}

Thus, the hyperbolic logarithm function $\ln$ is defined by the inverse isomorphism of $\mathrm{\exp}$.

\subsection{Bicomplex numbers} The algebra of bicomplex numbers defined in the introduction is the set
\begin{equation*}
\mathbb{B}: = \Big\{x+y\mathbf{i}+z\mathbf{j}+t \mathbf{k}:\, x,y,z,t \in \mathbb{R};\quad \ig, \jg,\kg \notin \R \Big\},
\end{equation*}
where $\ig, \jg, \kg$ are \textit{imaginary units} satisfying the following multiplication rules
\begin{equation*}
\mathbf{i}^{2} = \mathbf{k}^{2} = -1,\quad\mathbf{j}^{2} = 1,\quad\mathbf{i}\mathbf{k} = \mathbf{k}\mathbf{i} = \mathbf{j}.
\end{equation*}
 $\B$ contains three two-dimensional real subalgebras: two copies of the field of complex numbers, $\R(\ug): = \{x+ \ug y:\, x, y \in \R\}$, $(\ug = \ig, \kg)$ and the algebra of hyperbolic numbers $\R(\jg) = \D$. This implies that each bicomplex number $\omega$ has three $\R(\ug)$-algebraic representations, given by
 \begin{align}\label{u-representation}
 \omega = \Rj_\ug(\omega)+\pi(\ug) \Ij_\ug(\omega),
 \end{align}
 where $\pi$ is the permutation $\pi = \left(
\begin{array}
{ccc} \ig & \jg & \kg\\
 \kg& \ig & \jg
\end{array}
\right)$ and$\Rj_\ug(\omega), \Ij_\ug(\omega) \in \R(\ug)$ for each $\ug \in \{\ig,\jg,\kg\}$.

 Write $\omega = \Rj_\jg(\omega)+\ig \Ij_\jg(\omega)$ then, from the representation~\eqref{specdecomp} one can derive that
\begin{equation}\label{idemprepr}
 \omega = \mathcal{P}_{1}(\omega)\mathbf{e}_{1}+\mathcal{P}_{2}(\omega)\mathbf{e}_{2},
\end{equation}
where $\mathcal{P}_k$ is the algebra homomorphism form $\B$ to $\C = \R(\ig)$, defined by
\[
\mathcal{P}_k(\omega) = \pi_k(\Rj_\jg(\omega))+\ig \pi_k(\Ij_\jg(\omega)),\quad (k = 1,2).
\]
It follows that each bicomplex number can be viewed as a pair $(z_1,z_2)$ of the product algebra $\C\times \C$ via the map $\omega\mapsto (\mathcal{P}_1(\omega), \mathcal{P}_2(\omega))$. Thus the group $\B_\ast$ of units of $\B$ is the set of all bicomplex numbers $\omega$ such that $\mathcal{P}_{1}(\omega)\neq 0$ and $\mathcal{P}_{2}(\omega)\neq 0$.

On bicomplex numbers one has three conjugations $ \overline{\omega}^\ig$, $\overline{\omega}^\jg$ and $\overline{\omega}^\kg$ called the conjugations with respect to $\ig$, $\jg$ and $\kg$, respectively (\cite{DRA}). From the representation~\eqref{u-representation} we have
\begin{equation}\label{conjugations}
\overline{\omega}^\ug = \Rj_\ug(\omega)-\pi(\ug) \Ij_\ug(\omega), \quad(\ug = \ig, \jg, \kg).
\end{equation}
Therefore, from~\eqref{conjugations} and~\eqref{u-representation}, $ \omega \overline{\omega}^\ug$ belongs to $\R(\ug)$, $(\ug = \ig, \jg, \kg)$. In particular, $\omega \overline{\omega}^\jg \in \D^+ $. Thus, since $\D$ is square-root closed~\cite{HS} (i.e., every positive hyperbolic number $u$ has a unique positive square root $\sqrt{u}$) the $\jg$-modulus $|\omega|_\jg$ of $\omega$ is given by
\begin{align}\label{jmodulus}
        |\omega|_\jg: = \sqrt{\omega \overline{\omega}^{\jg}}.
\end{align}
Some remarkable properties of $|.|_\jg$ are given by the following statements. For the proof see~\cite[Section 4.4]{DRA}. Only for triangular inequality $M3)$ we refer to~\cite[Section 1.5]{SHAP}. For $\omega,\psi \in \B$,
 \begin{itemize}
 	\item [$M1)$] $|\omega|_\jg = 0 \mbox{~if and only if~} \omega = 0$;
 	\item [$M2)$] $ |\omega \psi|_\jg = |\omega|_\jg |\psi|_\jg$;
 	\item [$M3)$] $|\omega+\psi|_\jg\leq |\omega|_\jg+|\psi|_\jg$;
 	\item [$M4)$]$|\omega|_\jg = \|\mathcal{P}_{1}(\omega) \| \eg_1+\|\mathcal{P}_{2}(\omega) \| \eg_2$;
 	\item [$M5$)] $\|\omega\| = \sqrt{\Rj(|\omega|^2_\jg)}$,
 \end{itemize}
 where $\|.\|$ is the Euclidean norm on $\B$ which coincides with that in $\R(\ug)$, $(\ug = \ig, \jg, \kg)$ and with the modulus for $\ug = \ig, \kg$.

Finally, recall that a function $f:O\longrightarrow \B$ defined in the open set $O\subset \B$ is said to be $\B$-\ in $O$ if, for every $\omega \in O$ there exists a number $f'(\omega) \in \B$ such that
\begin{equation*}
f'(\omega): =	\lim_{\substack{\psi\mapsto \omega\\
	(\psi-\omega) \in \B_\ast}} \frac{f(\psi)-f(\omega)}{\psi-\omega}.
\end{equation*}
 For more details concerning bicomplex holomorphicity, we refer the reader to the following references: \cite{SHAP}, \cite{LUNA2} and \cite{SR}.

\section{Hyperbolic valued norm on bicomplex numbers}
The notion of \textit{$\D$-norm} on bicomplex numbers is introduced by Alpay et al. in~\cite{SHAP} and considered by Kumar et al.~\cite{C*alg} in the study of bicomplex $C^\ast$-algebra. In this section we establish additional properties for $\D$-norms on $\B$ viewed as the complexification of the $f$-algebra $\D$. Special attention is paid to the notion of Riesz subnorm of a $\D$-norm that plays a crucial role in the proof of Theorem~\ref{main}.

\subsection{Definition and properties}
Recall that a function $\Ng: \B \rightarrow \D$ is called a hyperbolic valued norm or $\D$-norm on $\B$ if the following properties are satisfied:
\begin{itemize}
	\item [\rm{(i)}] $\Ng(\omega) = 0$ implies $\omega = 0$;
	\item [\rm{(ii)}] $\Ng(\lambda \omega) = |\lambda| \Ng(\omega)$ for all $\lambda \in \mathbb{R}$ and for all $\omega \in \B$;
	\item [\rm{(iii)}] $\Ng(\omega+\psi) \leq \Ng(\omega)+ \Ng(\psi)$ for all $\omega, \psi \in \B$.
\end{itemize}

Clearly, every $\D$-norm $\Ng$ is positive, i.e., $\Ng(\omega) \in \D^+$ for all $\omega \in \B$. Moreover, one can see that the second triangular inequality holds.
\begin{proposition}
\begin{equation}\label{traigineq}
\Big|\Ng(\omega)-\Ng(\psi)\Big|\leq \Ng(\omega-\psi) \mbox{\quad for all~} \omega,\psi \in \B.
\end{equation}
\end{proposition}
\begin{proof}
Observing that $\pm(\Ng(\omega)-\Ng(\psi))\leq\Ng(\omega-\psi)$ one has, passing to supremum from the above,
$|\Ng(\omega)-\Ng(\psi)|\leq \Ng(\omega-\psi)\leq\Ng(\omega)+ \Ng(\psi).$
\end{proof}

For example, according to the properties $M1), M2)$ and $M3)$, the $\jg$-modulus defined in~\eqref{jmodulus} is a multiplicative $\D$-norm on $\B$ and satisfies the following properties.

\begin{proposition}\label{jhnorm} For all $\omega \in \B$ one has

\begin{itemize}
		\item[\rm ({i})]$ | \omega|_\jg = \sup \Big\{\Rj_\jg(\omega) \cos \theta + \Ij_\jg(\omega) \sin \theta:\, \theta \in \left[ 0, 2\pi\right] \Big\}$;
		\item[\rm ({ii})]$ |\omega|_\jg\geq |\Rj_\jg(\omega)|$ \quad and \quad $|\omega|_\jg\geq |\Ij_\jg(\omega)|.$

\end{itemize}
\end{proposition}
\begin{proof}
	(i): Follows from the closure of $\D$ for the square mean (\cite[Theorem 6.1]{HS}) by observing that $|\omega|_\jg = \sqrt{\Rj_\jg{(\omega)}^2 +\Ij_\jg{(\omega)}^2} $.

	(ii): The identity $|\omega|^2_\jg = \Rj_\jg{(\omega)}^2 + \Ij_\jg {(\omega)}^2$ implies that
\[
	|\omega|^2_\jg\geq \Rj_\jg{(\omega)}^2
	\quad \textup{and} \quad
	|\omega|^2_\jg\geq \Ij_\jg{(\omega)}^2.
\]
	So, since the hyperbolic square-root function is increasing and we have $\sqrt{z^2} = |z|$ for every $z\in \D$. Then,
	
\begin{equation*}
	|\omega|_\jg\geq |\Rj_\jg(\omega)| \mbox{~and }|\omega|_\jg\geq |\Ij_\jg(\omega)|.
\qedhere
\end{equation*}
\end{proof}
	 Recall that if $A+\ig A$ is the complexification of the unitary Archimedean $f$-algebra $A$ then for every $\mathfrak{z} = a+ \ig b \in A+\ig A$, the supremum
\begin{equation}\label{modulus}
			|\mathfrak{z}|: = \sup\{a \cos \theta+ b \sin \theta:\, \theta \in\left[ 0, 2 \pi\right] \}
\end{equation}
			exists in $A$ and is called the modulus of $\mathfrak{z}$ and satisfies the properties:
\begin{equation*}
			(i) |\mathfrak{z}| = 0 \mbox{~iff~} \mathfrak{z} = 0,\quad (ii) |\mathfrak{z} \mathfrak{w}| = |\mathfrak{z}| |\mathfrak{w}|, \quad (iii) |\mathfrak{z}|\wedge |\mathfrak{w}| = 0 \mbox{~iff~} \mathfrak{z} \mathfrak{w} = 0.
\end{equation*}
			$A+\ig A$ is said to be normal if $A+\ig A = \{\mathfrak{z}\}^{\perp}+ \{\mathfrak{w}\}^{\perp}$ for all $\mathfrak{z}, \mathfrak{w} \in A+ \ig A$ such that $ |\mathfrak{z}|\wedge | \mathfrak{w}| = 0$, where $\{\mathfrak{u}\}^{\perp} = \{\mathfrak{v}\in A+\ig A :\,|\mathfrak{u}|\wedge | \mathfrak{v}| = 0 \}$. In this case every $\mathfrak{z} \in A+\ig A$ has a polar-decomposition, i.e., there exists $\mathfrak{u} \in A+\ig A $ such that $\mathfrak{z} = \mathfrak{u} | \mathfrak{z}|$ and $|\mathfrak{z}| = \bar{\mathfrak{u}} \mathfrak{z}$; here $\bar{\mathfrak{u}}$ is the conjugate of $\mathfrak{u}$. More details about complexification of $f$-algebras can be found in~\cite{Beurk}.
			In the complexification $\B = \D+\ig \D$ the modulus~\eqref{modulus} of $\omega$ is its $\jg$-modulus $|\omega|_\jg$ (Proposition~\ref{jhnorm}).
			
\begin{proposition}\label{normaility}
				 $\B = \D+ \ig \D$ is normal.
			
\end{proposition}
		
\begin{proof}
		Let $\omega \in \B$ and let $ \psi,\varphi \in \B$ be such that $ | \psi|_\jg\wedge |\varphi|_\jg = 0 $, i.e., $\psi \varphi = 0$. If $\psi = 0$ we have $\{\psi\}^{\perp} = \B$ so for $\varphi_1 \in \{\varphi\}^{\perp}$ and $ \psi_1 \in \{\psi\}^{\perp}$ such that $\psi_1 = \omega-\varphi_1$ we have, $\omega = \varphi_1+\psi_1$. Similarly if $\varphi = 0$. Assume now that $\psi\neq0$ and $\varphi\neq0$. Therefore, $\varphi = z \eg$ and $\psi = w \bar{\eg}$ for some $z,w \in \C\setminus\{0\}$ and some $\eg = \eg_1, \eg_2$. Thus, $ \{\varphi\}^{\perp} = \{z \bar{\eg}:\, z \in \C\}$ and $\{\psi\}^{\perp} = \{z\eg:\, z\in \C\}$. Which implies from the idempotent representation~\eqref{idemprepr} that $\omega = \varphi_1+ \psi_1$ with $\psi_1\in \{\psi\}^{\perp}$ and $\varphi_1 \in \{\varphi\}^{\perp}$.
\end{proof}

		As mentioned above, normality yields polar decomposition. The following result is a direct consequence of Proposition~\ref{normaility}.
		
\begin{proposition}[Polar decomposition ]\label{polar-decomp}
			For every $\omega \in \B$ there exist $\upsilon \in \B$ such that
\begin{equation*}
			\omega = \upsilon |\omega|_\jg \mbox{~} (\mbox{and~} |\omega|_\jg = \overline{\upsilon}^{\jg} \omega).
\end{equation*}
\end{proposition}

\begin{remark}\label{remquiv} \rm {Obviously, the function $\Ng = N_1 \eg_1+N_2 \eg_2$ is a $\D$-norm on $\B$ whenever $N_1$ and $N_2$ are two real norms on $\B$. But the converse is false. The components of a $\D$-norm (in the spectral decomposition) are semi-norms. Indeed, from $M4)$ the components of the $\D$-norm $|.|_\jg$ are $\|\mathcal{P}_1(.)\|$ and $\|\mathcal{P}_2(.)\|$. However, $\|\mathcal{P}_1(\eg_2)\| = \|\mathcal{P}_2(\eg_2)\| = 0 $.}
\end{remark}
The above remark leads to the following definition
\begin{definition}
	A $\D$-norm $\Ng = N_1 \eg_1+N_2 \eg_2$ on $\B$ is said to be integral if the semi-norms $N_1$ and $N_2$ are norms on
	 $\B$.
\end{definition}
 From the above definition a $\D$-norm $\Ng$ on $\B$ is integral if and only if $\|\Ng(\omega)\|_h = 0$ implies $\omega = 0$.

 Now, equivalence of $\D$-norms on $\B$ is defined as
\begin{definition}\label{equivnorms}
	Let $G$ be the group constituted by the identity and the conjugation operator of $\D$. Two $\D$-norms $\Ng_1$ and $\Ng_2$ on $\B$ are called equivalents if there exists two real numbers $k, k'> 0$ and $\mathcal{L} \in G$ such that
	
\begin{equation*}
	k' \mathcal{L} (\Ng_2(\omega)) \leq \Ng_1(\omega) \leq k \mathcal{L}(\Ng_2(\omega)) \quad \mbox{for all} \omega \in \B.
\end{equation*}
		
\end{definition}
Let us introduce the binary relation $\sim$ defined for a pair of $\D$-norms by
\begin{equation*}
\Ng_1\sim \Ng_2 \mbox{~if and only if~} \Ng_1 \mbox{~and~} \Ng_2 \mbox{~are equivalent}.
\end{equation*}
It is easy to verify that $\sim$ is an equivalence relation. Note that $\Ng_1\nsim \Ng_2 $ when $\Ng_1$ is integral and $\Ng_2$ is not. Indeed, suppose the contrary, then there exist a real $k> 0$ and $ \mathcal{L}\in G$ such that $ \Ng_1(\omega) \leq k \mathcal{L}(\Ng_2(\omega))$ for all $\omega \in \B$. Since $\Ng_2$ is not integral then there exist a nonzero $\omega_0 \in \B$ such that $ \|\Ng_2(\omega_0)\|_h = \| k \mathcal{L}(\Ng_2(\omega_0))\|_h = 0$. This implies that $ k \mathcal{L}(\Ng_2(\omega_0))
\in \eg \R $ for some $\eg = \eg_1,\eg_2$. It follows that $\Ng_1(\omega_0) \in \eg \R$, since $\eg \R $ is an order ideal in $\D$ (see~\cite{HS}). Therefore, $\|\Ng_1(\omega_0)\|_h = 0 $. Which is a contradiction, since $\Ng_1$ is integral.
\begin{proposition}
	Two integral $\D$-norms on $\B$ are equivalent.
\end{proposition}
\begin{proof}
	Straightforward since real norms are equivalent in finite dimensional vector spaces.
\end{proof}

In the following proposition we introduce our main tool in the study of the convergence of the bicomplex zeta and gamma function.

\begin{proposition}[Riesz subnorm]\label{N-norm}

Let $\Ng$ be a $\D$-norm on $\B$. Then the function $\lceil \Ng \rceil:\B\longrightarrow \R$ defined by
\begin{equation*}
\lceil \Ng \rceil(\omega): = \min \{\alpha \in \R^+ :\, \alpha\geq \Ng(\omega)\} = \Ng(\omega) \vee \overline{\Ng(\omega)}
\end{equation*}
is a real norm on $\B$, called Riesz subnorm of $\Ng$, that satisfies the following properties for all $\varphi, \psi \in \B:$
\begin{itemize}
	\item [$(\rm {i})$]$\lceil \Ng \rceil(\varphi) \leq \lceil \Ng \rceil(\psi) \quad \text{whenever} \quad \Ng(\varphi)\leq \Ng(\psi)$;
	\item [$(\rm {ii})$]$ \Ng(\varphi)\leq \lceil \Ng \rceil(\varphi)$.
\end{itemize}
\end{proposition}
\begin{proof}
 The proof follows immediately from the properties of the norm $\|.\|_R$~\eqref{riesznorm}, by observing that $\lceil \Ng \rceil(\omega) = \| \Ng(\omega)\|_R$ for all $\omega \in \B$.
\end{proof}
Applying $\|.\|_R$ in triangle inequality~\eqref{traigineq} one can derive that every $\D$-norm $\Ng$ is a Lipschitz function from $(\B, d_{\lceil \Ng \rceil}) $ to $(\D, d_R)$ where $d_{\lceil \Ng \rceil}$ and $d_R$ are the metrics defined by the norms $\lceil \Ng \rceil $ and $\|.\|_R$ respectively.

Throughout the paper we will write
\begin{equation}\label{SubN}
\lceil \Ng \rceil(\omega) = \|\omega\|_\jg \quad \mbox{for the~} \mbox{$\D$-norm}\quad \Ng(\omega) = |\omega|_\jg.
\end{equation}

\begin{proposition}\label{jnorme}
	For every $\omega,\psi \in \B$

\begin{itemize}
		\item [$(\rm {i})$] $\|\omega \|\leq\|\omega\|_\jg \leq \sqrt{2}\|\omega\| $;
		\item [$(\rm {ii})$] $\|\omega \psi\|_\jg\leq \|\omega\|_\jg \|\psi\|_\jg$, with equality whenever $\omega \in \R(\ug)$ or $\psi \in \R(\ug)$, $(\ug = \ig, \kg)$. Thus, the real norm $\|.\|_\jg$ is submultiplicative.
\end{itemize}
\end{proposition}

\begin{proof}
For the proof we will use the following two elementary properties that hold in any $f$-algebra. For $u,v\geq 0$,
	
\begin{itemize}
		\item [	$P1)$] 	 ${(u\vee v)}^2 = u^2 \vee v^2;$
		\item [	$P2)$] $u\vee v\leq u+v.$
	
\end{itemize}
$(\rm {i})$ Let $\omega\in \B$. By definition of $\|\omega\|_\jg$ (Proposition~\ref{N-norm}), property $P1)$ implies that
\begin{equation}\label{eqmodulus1}
 \|\omega\|_\jg^2 = |\omega|_\jg^2 \vee \overline{|\omega|_\jg}^2.
\end{equation}
 So, since the mapping $z\mapsto \Rj(z)$ is a positive operator from $\D$ to $\R$ (i.e., a linear form such that $\Rj(z)\leq \Rj(w)$ for all $z,w \in \D$ with $z\leq w$) then from Eq\eqref{eqmodulus1} and $M5)$ (Section 2) one has
\begin{equation*}
\|\omega\|_\jg^2 \geq \Rj(|\omega|_\jg^2) = \|\omega\|^2.
\end{equation*}
For the second inequality, Eq~\eqref{eqmodulus1} together with $P2)$ yields that $\|\omega\|^2_\jg \leq |\omega|_\jg^2 +\overline{|\omega|_\jg}^2$ and then $\|\omega\|_\jg^2 \leq 2\Rj(|\omega|_\jg^2)$. Hence, again by $M5)$ we obtain
\begin{equation*}
\|\omega\|_\jg^2 \leq 2\|\omega\|^2.
\end{equation*}
$(\rm {ii})$ Let $\psi,\omega \in \B$. It follows by $N2)$ (section 2) and by $\| \omega \psi\|_\jg = \| |\omega \psi|_\jg\|_R$ that we have
\begin{equation*}
	\|\omega \psi\|_\jg \leq \|\omega\|_\jg \| \psi\|_\jg.
\end{equation*}
	Suppose now that $\omega \in \R(\ug)$, $(\ug = \ig,\kg)$. Then $\omega = a+ \varepsilon_{\ug} \ig b$ where $a,b \in \R$ and $\varepsilon_{\ug}$ be such that $\varepsilon_{\ig} = 1$ and $\varepsilon_{\kg} = \jg$. Then $|\omega|_\jg = \sqrt{a^2 +b^2} = \|\omega\|$. Which implies that $\|\omega \psi\|_\jg = \|\omega\|_\jg \|\psi\|_\jg$.
\end{proof}
\begin{remark}
	\rm {Let $\mathrm{S}$ and $\mathrm{S}_\jg$ be the unit spheres in $(\B,\|.\|)$ and $(\B, \|.\|_\jg)$, respectively. We have that $\|\sqrt{2} \eg_1\|_\jg = \sqrt{2}$ and $\|1\| = 1$ with $\sqrt{2} \eg_1 \in \mathrm{S}$ and $1 \in\mathrm{S}_\jg $. This implies by inequality $(\rm {i})$ (Proposition~\eqref{jnorme}) that
		
\begin{equation*}
		\sqrt{2} = \displaystyle \sup_{\omega \in \mathrm{S}} \|\omega\|_\jg \mbox{~and~} 1 = \displaystyle \sup_{\omega \in \mathrm{S}_\jg} \|\omega\|.
\end{equation*}
}

\end{remark}
Property (\rm{ii}) of proposition~\ref{jnorme} means that $(\B,+,., \|.\|_\jg)$ is a real Banach algebra. Thus, for every bicomplex number $\omega$ the exponential of $\omega$ can be defined as the absolute convergence series given by
\begin{equation}\label{seriexp}
e^\omega: = \displaystyle\displaystyle\sum_{n = 0}^\infty \frac{\omega^n}{n!} = e^{\mathcal{P}_1(\omega)}\eg_1+ e^{\mathcal{P}_2(\omega)}\eg_2.
\end{equation}
The bicomplex exponential function $\mathrm{Exp}$ is a group homomorphism from the additive group $\B$ to the multiplicative group $\B_\ast$ with $\mathrm{ker(Exp)} = 2\ig \pi \ZG$, where
\begin{equation}\label{hypint}
\ZG: = \Z \eg_1+ \Z \eg_2.
\end{equation}
$\ZG$ is a sublattice and subring of $\D$, called the ring of hyperbolic integers~\cite{HS2}.
\subsection{Topology of Bicomplex numbers}
Let $\Ng$ be a $\D$-norm on $\B$ and let $\omega_0 \in \B$ and $ \tau \gg 0$. Define the open $\D$-ball $\Bg_\Ng^o (\omega_0,\tau) $ and the closed $\D$-ball $\Bg_\Ng(\omega_0,\tau)$ centered at $\omega_0$ with hyperbolic radius $\tau$ by
\begin{align*}
\Bg_\Ng^o (\omega_0,\tau): &=& \Big\{\omega \in \B:\, \Ng(\omega-\omega_0) \ll \tau\Big\};\\
\Bg_\Ng(\omega_0,\tau): &=& \Big\{\omega \in \B:\, \Ng(\omega-\omega_0) \leq \tau\Big\}.
\end{align*}
 If $\Ng$ is a real function (i.e., $\Ng(\omega) = \overline{\Ng(w)}$ for all $\omega \in \B$) so that $\Ng$ is a real norm then $\Bg_\Ng^o (\omega_0,\tau))$ and $\Bg_\Ng(\omega_0,\tau) $ are usual balls in the normed space $(\B, \Ng)$ with radius $r = \tau\wedge \bar{\tau} \in \R^+_\ast$.
\begin{proposition}\label{top}
	Every $\D$-norm $\Ng$ on $\B$ generate a topology $\mathcal{T}_{\Ng}$ defined as follows: a nonempty subset $\mathcal{O}$ of $\B$ is said to be open if, each point of $\mathcal{O}$ is the center of some $\D$-ball.
\end{proposition}
\begin{proof}
	For the proof, it suffice to shows that the intersection of finite open sets is an open set. This follows from the closure of $\D_\ast^+$ under the lattice operation $\wedge$ by observing that $\displaystyle \bigcap_{i = 1}^{n}\Bg_\Ng(\omega,\tau_i) = \Bg_\Ng(\omega,\displaystyle\bigwedge_{i = 1}^{n} \tau_i)$.
\end{proof}
From this topology the limit in bicomplex numbers can be formulated by
\begin{proposition}
	Let $f: (\B,\mathcal{T}_{\Ng_1})\longrightarrow (\B,\mathcal{T}_{\Ng_2})$. Then, a bicomplex number $\psi $ is the limit of $f$ at the point $\omega_0 \in \B$ if and only if for every hyperbolic number $\xi\gg 0$ there exists a hyperbolic number $ \eta \gg 0$ such that
\begin{equation*}
	\Ng_2(f(\omega)- \psi)\leq \xi \quad \mbox{whenever} \quad \Ng_1(\omega- \omega_0)\leq \eta.
\end{equation*}
\end{proposition}

\begin{proposition}\label{topcanonic} $\D$-norms on $\B$ are
topologically equivalent.
\end{proposition}
\begin{proof}
	Let $\Ng$ be a $\D$-norm on $\B$. Since all real norms on finite dimensional real vector space are topologically equivalent it suffice to prove that $\mathcal{T}_\Ng$ is equivalent to $\mathcal{T}_{\lceil \Ng \rceil}$. One has to prove that $Id:(\B,\mathcal{T}_\Ng) \longrightarrow (\B,\mathcal{T}_{\lceil \Ng \rceil})$ and its inverse $Id^{-1}:(\B,\mathcal{T}_{\lceil \Ng \rceil}) \longrightarrow (\B,\mathcal{T}_{\Ng})$ are continuous. Let $\psi,\varphi \in \B$. For $\epsilon \in \R^+_\ast$ take $\eta = \epsilon$ we have
\begin{equation*}
	\Ng(\psi-\varphi)\leq \eta \mbox{~implies~}\lceil \Ng \rceil(\psi-\varphi) = \|\Ng(\psi-\varphi)\|_R\leq\epsilon.
\end{equation*}
	Conversely, For $ \xi\gg 0$, put $ \eta = \xi \wedge \bar{\xi} \in \R^+_\ast$. Thus one has
\begin{equation*}
	\lceil \Ng \rceil(\omega-\psi)\leq \eta \mbox{~implies~} \Ng(\omega-\psi)\leq\xi.
\qedhere
\end{equation*}
\end{proof}

\begin{remark}
We know from Remark~\ref{remquiv} that two $\D$-norms on $\B$ are not necessarily equivalent in sens of Definition~\ref{equivnorms}. Nevertheless (by Proposition~\ref{topcanonic}) they are topologically equivalent. Moreover, one can see that $\D$-bounded sets in a $\D$-normed space $(\B, \Ng)$ are bounded sets in the real normed space $\B$.
\end{remark}

\section{Bicomplex trigonometry}
In this section we develop the concept of $\D$-trigonometric form of a nonzero bicomplex number introduced by Luna-Elizarrar\'as et al.~\cite[Chapitre 3]{LUNA1}. Using the $f$-algebra structure of $\D$ one is able to select a specified $\D$-valued modulus and argument of bicomplex numbers.

\subsection{Basic concepts and properties}
As in complex numbers, bicomplex trigonometry in his basic form is the study of the properties of hyperbolic cosine and sine functions. From the Banach algebra structure of $\D$, circular functions can be defined for all $z \in \D$ as
\begin{align*}
 \cos (z)&: = \sum_0^\infty \frac{{(-1)}^n z^{2 n}}{2 n!} = \cos(\pi_1(z)) \eg_1+ \cos (\pi_2(z)) \eg_2,\\
 \sin (z)&: = \sum_0^\infty \frac{{(-1)}^n z^{2 n+1}}{(2 n+1)!} = \sin(\pi_1(z)) \eg_1+ \sin (\pi_2(z))\eg_2.
\end{align*}
We will give some of the properties of cosine and sine. The proof follows immediately from usual properties of the real cosine and sine functions, using the above spectral decompositions. For $z,w \in \D$; $\varepsilon \in \SG$; $h \in \ZG$,
\begin{itemize}
	\item [$C1$)] $e^{iz} = \cos z+ \ig \sin z$.
	\item [$C2)$] $\cos (z+2 \pi h) = \cos (z)$ and $\cos (z+2 \pi h) = \cos (z) $.
	\item [$C3)$] $\cos (\varepsilon z) = \cos(z)$ and $\sin (\varepsilon z) = \varepsilon\sin(z) $.
	\item [$C4)$] $z^2 +w^2 = 1$ if and only if $z = \cos \theta$ and $w = \sin \theta$ for some $\theta \in \D$.
		\item [$C5)$] The restriction of the cosine function to $[ 0,\pi]_\D$ establish a bijection with $[-1,1]_\D$, its inverse is denoted $\arccos$.
	\item [$C6)$] The restriction of the sine function to $[ \frac{-\pi}{2},\frac{\pi}{2}]_\D$ establish a bijection with $[-1,1]_\D$, its inverse is denoted $\arcsin$.
\end{itemize}

\begin{proposition}\label{thtrigo}
Every nonzero bicomplex number $ \omega$ can be written in the form
\begin{equation*}
	\omega = |\omega|_{\mathbf{j}}\left(\cos\phi+\mathbf{i}\sin\phi\right),
\end{equation*}
where $\phi$ is an hyperbolic number called a $\D$-argument of $\omega$.
\end{proposition}
\begin{proof}
Let $\omega$ be a nonzero bicomplex number. If $\omega$ is invertible then from polar decomposition (Proposition~\ref{polar-decomp}), $\omega$ can be written as
\begin{equation*}
\omega = |\omega|_\jg \upsilon,
\end{equation*}
where $\upsilon \in \B$ with $|\upsilon|_\jg = 1$ that is, $\Rj_\jg{(\upsilon)}^2 +\Ij_\jg{(\upsilon)}^2 = 1$ which implies from property $C5)$ above that $\upsilon = \cos \phi+ \ig \sin \phi$ for some $\phi\in \D$. If $\omega$ is non-invertible, i.e., $\omega = z \eg$ for some nonzero complex number $z$ and some $\eg \in\{\eg_1,\eg_2\}$, then there exists $\phi \in \R\subset \D $ such that
\begin{equation*}
\omega = \|z\| (\cos \phi+ \ig \sin \phi)\eg = |\omega|_\jg (\cos \phi+ \ig\sin \phi).
\qedhere
\end{equation*}
\end{proof}

\begin{proposition}\label{PArg} Let $\omega \in \B_\ast$ then the set $\ar_\D(\omega)$ of all $\D$-arguments of $\omega$ has a unique element $\phi_p \in \left(-\pi,\pi\right] _\D$ called principal $\D$-argument of $\omega$, denoted $\Ar_\D(\omega)$.
\end{proposition}
\begin{proof}
Let $\omega \in \B_\ast$. So by Proposition~\ref{thtrigo} and $C1)$, $ \phi_1,\phi_2 \in \ar_\D(\omega)$ if and only if $e^{\ig \phi_1} = e^{\ig \phi_2}$, i.e., if and only if $\ig(\phi_1-\phi_2) \in \mathrm{ker}(\mathrm{Exp}) = \ig 2 \pi \ZG$. Hence
\begin{equation*}
\ar_\D(\omega) = \phi_0+2 \pi \ZG, \mbox{\quad for some~}\phi_0 \in \ar_\D(\omega).
\end{equation*}
It follows that $\ar_\D(\omega)$ has a unique element $\phi_p$ satisfying $-\pi \ll \phi_p\leq \pi$.
\end{proof}
The principal $\D$-argument of $\omega\in \B_\ast$ is determined by the solution of the equations
\begin{equation}\label{eqcos}
\cos \phi = \frac{\Rj_\jg(\omega)}{|\omega|_\jg} \mbox{~and~} \sin \phi = \frac{\Ij_\jg(\omega)}{|\omega|_\jg},\quad \phi \in \left(-\pi,\pi \right] _\D.
\end{equation}
For example, $\Ar_\D(1) = 0, \Ar_\D(\ig) = \frac{\pi}{2}, \Ar_\D(\jg) = \pi \eg_1$ and $\Ar_\D(\kg) = -\frac{\pi}{2} \jg$. Notice that the existence of the principal $\D$-argument is not guaranteed when $\omega$ is a zero divisor. For example, the two $\D$-arguments $0$ and $\pi \eg_1$ of $\eg_1$ belongs to $\left(-\pi,\pi \right] _\D.$

The following properties are direct consequence of proposition~\ref{PArg} (equalities are modulo $2\pi \ZG$).\newpage

\begin{corollary}
For $\varphi,\psi \in \B_\ast$ one has :
\begin{itemize}
	\item [$A1)$] $\Ar_\D(\overline{\varphi}^\ig) = \overline{\Ar_\D(\varphi)}$, \quad $\Ar_\D(\overline{\varphi}^\kg) = -\overline{\Ar_\D(\varphi)} $, \quad $ \Ar_\D(\overline{\varphi}^\jg) = -\Ar_\D(\varphi).$
	\item [$A2)$] $\Ar_\D(\varphi^{-1}) = -\Ar_\D(\varphi)$, \quad $\Ar_\D(\varphi \psi) = \Ar_\D(\varphi)+\Ar_\D(\psi)$.
\end{itemize}
\end{corollary}

One can now give the following definition.
\begin{definition}
The principal branch of the bicomplex Logarithm of $\omega \in \B_\ast$ is given by
\begin{equation*}\label{logbranch}
 \mathrm{Log}(\omega): = \ln (|\omega|_\jg)+ \ig\Ar_\D(\omega),
\end{equation*}
and for $\alpha \in \B$, the bicomplex exponentiation
\begin{equation}\label{bicompawer}
\mathrm{Exp}_\omega(\alpha) = \omega^\alpha: = e^{\alpha \mathrm{Log}(\omega)}.
\end{equation}
\end{definition}
\begin{proposition}
 For $\varphi,\psi,\omega \in \B_\ast, \alpha \in \B$ and $\ug = \ig, \kg, \jg$, one has (equalities $L1$ and $L2$ are modulo $2\ig \pi \ZG$)
\begin{itemize}
	\item [$L1)$]$\mathrm{Log}(\overline{\varphi}^{\ug}) = \overline{\mathrm{Log}(\varphi)}^{\ug}.$
	\item [$L2)$]$\mathrm{Log}(\varphi^{-1}) = -\mathrm{Log}(\varphi)$,\, $\mathrm{Log}(\varphi \psi) = \mathrm{Log}(\varphi)+\mathrm{Log}(\psi)$.
	\item
	[$L3)$]$\omega^{-\alpha} = \frac{1}{\omega^\alpha}$, $\omega^\alpha \omega^\beta = \omega^{\alpha+\beta}$ and ${(\varphi \psi)}^\alpha = \varphi^\alpha\psi^\alpha $ if and only if $\alpha \in \ZG$.
\end{itemize}
\end{proposition}
\begin{proof}
$L3)$ is straightforward.
$L1)$ and $L2)$ follows from $A1)$ and $A2)$.
\end{proof}

In order to give a geometrical interpretation of the principal $\D$-argument, we introduce the function $<.,. >_\jg: \B\times \B\longrightarrow \D$ defined by $<\varphi,\psi >_\jg: = \Rj_\jg \left(\varphi\, \overline{\psi}^\jg\right)$. The inner product $<.,.>_\jg$ is symmetric, $\D$-bilinear, positive-definite, i.e., $<\omega,\omega >_\jg\geq 0$ with equality if and only if $\omega = 0$, and satisfies the hyperbolic Cauchy-Schwartz inequality:
\begin{equation}\label{Chauchyineq}
|<\varphi,\psi >_\jg| \leq |\varphi|_\jg |\psi|_\jg \mbox{\quad for all~} \varphi,\psi \in \B.
\end{equation}
Thus from $C5)$ and~\eqref{Chauchyineq} we can define the $\D$-angle between two invertible bicomplex numbers $\psi$ and $\varphi$, by the formula
\begin{equation*}
		\mathrm{angl}_\D(\varphi,\psi): = \arccos \Big(\frac{<\varphi,\psi >_\jg}{|\varphi|_\jg |\psi|_\jg}\Big).
\end{equation*}
 Let $\omega \in \B_\ast$. By (Theorem of signs~\cite{HS}) we can write $\Ar_\D(\omega) = \varepsilon |\Ar_\D(\omega)|$ for some $\varepsilon \in \SG$. So, by $C3)$ $\cos(|\Ar_\D(\omega)|) = \cos(\Ar_\D(\omega))$ with $|\Ar_\D(\omega)| \in \left[ 0,\pi\right] _\D$ which implies from Eq~\eqref{eqcos} and the above formula that
\begin{equation*}
|\Ar_\D(\omega)| = \mathrm{angl}_\D(\omega,1).
\end{equation*}
For example, $\mathrm{angl}_\D(\ig,1) = \mathrm{angl}_\D(\kg,1) = \frac{\pi}{2}$ and $\mathrm{angl}_\D(\jg,1) = \pi \eg_1 $.

\begin{proposition}
	Let $ \varepsilon \in \SG$ and $\omega \in \B_\ast$ then $\omega \in \D^\varepsilon \mbox{~if and only if~}\Ar_\D(\omega) = \left(\frac{1-\varepsilon}{2}\right) \pi $.
\end{proposition}
\begin{proof}
Let $\varepsilon \in \SG$ and $\omega \in \B_\ast$. Since $\phi_\varepsilon = :\left(\frac{1-\varepsilon}{2}\right) \pi$ satisfies $-\pi \ll \phi_\varepsilon\leq \pi$, then from $C3)$ and Eq~\eqref{eqcos}, $\omega \in\D^\varepsilon$ if and only if $\cos\phi_\varepsilon = \frac{\Rj_\jg(\omega)}{|\omega|_\jg} \mbox{~and~} \sin \phi_\varepsilon = \frac{\Ij_\jg(\omega)}{|\omega|_\jg}$, i.e. if and only if $\Ar_\D(\omega) = \phi_\varepsilon$.
\end{proof}

\subsection{The $n^{th}$ roots of a bicomplex number.}
Let $n\geq 1$. Then each bicomplex number $\omega = z_1 \eg_1+z_2 \eg_2$, $(z_1,z_2 \in \C)$ has $n^{(2-\nu(\omega))}$ $n^{th}$ roots, where $\nu(\omega)$ is the number of $k\in \{1,2\}$ such that $z_k = 0 $. If $\omega$ is invertible, i.e. $z_k\neq0$, $(k = 1,2)$ then, as in~\cite[Section 6.4]{LUNA1}, the $n^2$ $n^{th}$ roots of $\omega$ are the numbers
\begin{equation*}
\omega_{(h_1,h_2)} = \sqrt[n]{|z_1|} e^{\ig (\frac{2 \pi h_1+ \Ar(z_1))}{n}} \eg_1+ \sqrt[n]{|z_2|} e^{\ig (\frac{2 \pi h_2+ \Ar(z_2))}{n}} \eg_2; \mbox{\quad }h_1,h_2 = 0, \cdots,n-1.
\end{equation*}

In the following, we give additional properties of the $n^{th}$ roots of an invertible bicomplex number.
Using hyperbolic integers $\ZG$ (\cite{HS2}), the $n^2$ $n^{th}$ roots are given by
\begin{equation*}
\displaystyle \omega_h = \sqrt[n]{|\omega|_\jg}\displaystyle \, e^{\ig(\frac{2\pi h+\Ar_\D(\omega)}{n})}:\quad h \in \ZG, \mbox{~} 0\leq h\leq n-1.
\end{equation*}
In particular, the $n^{th}$-roots of unity are
described by the following proposition.

\begin{proposition}\label{group}
\begin{itemize}
		\item [$(\rm {i})$]
		The $n^{th}$-roots of unity are the set
\begin{equation*}
		 \mathcal{U}_n: = \Big\{\upsilon_h = e^{\frac{2 \ig \pi h}{n}}: \quad h \in \ZG, \mbox{~}0 \leq h\leq n-1\Big\}.
\end{equation*}
		It is a subgroup of the unit $\D$-sphere $\mathrm{S_\D}: = \{\omega \in \B:\, |\omega|_\jg = 1\}$ of the $\D$-normed space $(\B,|.|_\jg)$.
		\item [$(\rm {ii})$] $\|\upsilon_p-\upsilon_q\| = 2 \| \sin\frac{(p-q)\pi}{n}\|$ for all $\upsilon_p,\upsilon_q \in \mathcal{U}_n$;
		\item [$(\rm {iii})$] $\displaystyle\displaystyle \sum_{\upsilon \in \mathcal{U}_n} \upsilon = 0 $ and $\displaystyle\displaystyle \prod_{\upsilon \in \mathcal{U}_n} \upsilon = 1$.
\end{itemize}
\end{proposition}

\begin{proof}
	(i): Straightforward.

	(ii): Follows from $M5)$ using the identity: $e^{\ig z}-e^{\ig w} = 2 \ig e^{\ig \frac{(z+w)}{2}} \sin \frac{(z-w)}{2}$ holds for all $z,w \in \D$.

    (iii): For the sum one has
\begin{align*}
		\sum_{\upsilon \in \mathcal{U}_n} \upsilon 
		&= \sum_{\substack{h \in \ZG \\ 0\leq h\leq n-1}}e^{\frac{2 \ig \pi h}{n}}\\
		&= n \sum_{0\leq h_1\leq n-1}e^{\frac{2 \ig \pi h_1}{n}} \eg_1+n
		\sum_{0 \leq h_2\leq n-1}e^{\frac{2 \ig \pi h_2}{n}} \eg_2 \\
		&= 0 \eg_1+ 0 \eg_2 = 0.
\end{align*}
	For the product, we have
\[
\prod_{\upsilon \in \mathcal{U}_n} \upsilon = e^{\frac{2 \ig \pi}{n} \sigma_n}, \qquad
\sigma_n = \sum_{\substack{h \in \ZG\\ 0 \leq h\leq n-1}} h.
\]
	We have
\begin{align*}
		\sigma_n
		&= n \sum_{0\leq h_1 \leq n-1} h_1 \eg_1+ n\sum_{0 \leq h_2\leq n-1} h_2 \eg_2 \\
		&= \displaystyle \frac{n^2 (n-1)}{2} \eg_1 +\frac{n^2 (n-1)}{2} \eg_ 2 \\
		&= \frac{n^2 (n-1)}{2}.
\end{align*}
	Hence
\begin{equation*}
	\displaystyle\displaystyle \prod_{u \in \mathcal{U}_n} u = e^{n (n-1) \ig \pi} = 1.
\qedhere
\end{equation*}
\end{proof}

In the complex plane the $n^{th}$ roots of unity are the vertices of a regular polygon inscribed in the euclidean circle $S^1 \backsimeq\R/\Z$. For bicomplex numbers the $n^{th}$ roots of unity are in the unit $\D$-sphere $\mathrm{S_\D}$ which has the following topological structure.
\begin{proposition}\label{torus}
	The unit $\D$-sphere $\mathrm{S_\D} = \{\omega \in \B:\, |\omega|_\jg = 1\}$ is homeomorphic to the two dimensional torus $\mathrm{T}^2 = \mathrm{S}^1 \times \mathrm{S}^1$.
\end{proposition}
\begin{proof}
  It follows from $\D/\ZG\backsimeq \R/\Z\times \R/\Z$ that $\mathrm{S}_\D$ is the two-dimensional torus $\mathrm{T}^2 = \mathrm{S}^1 \times \mathrm{S}^1$ via the homeomorphism $\bar{f}:\D/\ZG\longrightarrow \mathrm{S}_\D:\hat{z}\mapsto f(z)$ where $f$ is the continuous group homomorphism from $\D$ to $\B_\ast$ defined by $f(z) = e^{2 \pi \ig z}$, with $\mathrm{ker}(f) = \ZG$, $\mathrm{Im}(f) = \mathrm{S}_\D$ and ${(\bar{f})}^{-1}(\omega) = \widehat{\frac{1}{2\ig \pi} \mathrm{Log}(\omega)}$ for all $\omega \in \mathrm{S}_\D$.
\end{proof}

Recall that a \textit{toroid} is an ordinary polyhedron, topologically torus-like. Its Euler number is then $v-e+f = 0$. A toroid is said to be \textit{regular} if the same number of edges meet at each
vertex, and each face has the same number of edges. A toroid is \textit{in class $T_2$} if each face has four edges and at each vertex exactly four edges meet. For regular toroids and their classification see~\cite{SZ} and references therein.

Let $n\geq 3$. In view of Proposition~\ref{torus} the bicomplex $n^{th}$-roots of unity can be identified with the vertices of a regular toroid in class $T_2$ that has $n^2$ quadrilateral faces. In particular, a regular toroid with minimal faces in class $T_2$ has $3\times 3$ quadrilateral faces. It corresponds to the group $\mathcal{U}_3$ illustrated by Figure 1(a).

\begin{figure}[!ht]
	\centering
	\includegraphics[scale = 0.3]{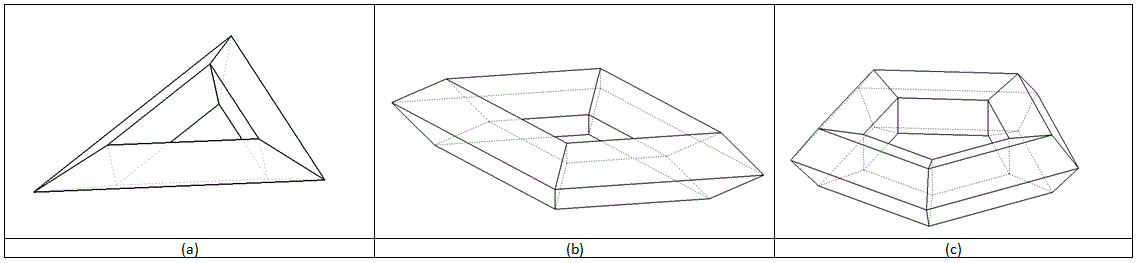}
	\caption{Toroids (a), (b) and (c) represent in class $T_2$ the groups $\mathcal{U}_3$, $\mathcal{U}_4$ and $\mathcal{U}_5$, respectively.}
\end{figure}

\newpage

\section{Proof of Theorem~\ref{main}}
The main purpose of this section is to prove Theorem~\ref{main}. We start by the absolute and uniform convergence of bicomplex zeta function. We introduce the bicomplex gamma function as an absolute convergent integral.
The Mellin transform of a bicomplex-valued function was first considered in~\cite{AGA} as a pair of complex Mellin transform using idempotent representation. In this section we extend Mellin integral to bicomplex numbers using absolute convergence of the integral. To the best of our knowledge this result and the method have never been considered before.

First of all we need some preliminary notations and properties.
Let us introduce two functions defined on $\B$ with values in $\R$ by
\begin{equation}
s_\R(\omega) = \Rj_\jg (\omega)\vee\overline{\Rj_\jg (\omega)}\quad \mbox{and} \quad i_\R(\omega) = \Rj_\jg (\omega)\wedge \overline{\Rj_\jg (\omega)}.
\end{equation}

\begin{proposition}
For $\omega \in \B$ and $\lambda \in \R$,
\begin{itemize}

	\item [$R1)$]
$s_\R$ and $i_\R$ are continuous and surjective.
	\item [$R2)$] $s_\R(\lambda \omega) = \lambda^{+}s_\R(\omega)-\lambda^{-} i_\R(\omega) \mbox{~and~}i_\R(\lambda \omega) = \lambda^{+}i_\R(\omega)-\lambda^{-} s_\R(\omega)$.
	\item [$R3)$] $i_\R(\omega+\lambda) = i_\R(\omega)+\lambda$ and $s_\R(\omega+\lambda) = s_\R(\omega)+\lambda $
	\item [$R4)$] $i_\R(\omega)> \lambda$ in $\R$ if and only if $\Rj_\jg(\omega)\gg \lambda$ in $\D$.

\end{itemize}
\end{proposition}
\begin{proof}
	$R1) $ The continuity follows from the continuity of $\Rj_\jg (.)$. The surjectivity holds since $s_\R(x) = i_\R(x) = x$ for all $x \in \R$.\\
	$R2) $ Follows from the identities: $u z \vee u w = u^+ (z\vee w)-u^{-}(z\wedge w)$ and $u z \wedge u w = u^+ (z\wedge w)-u^{-}(z\wedge w) $ that hold in $\D$ (see~\cite{HS}).\\
	$R3)$ Deduced from the identities: $u+(z\vee w) = (u+z)\vee(u+w)$ and $u+(z\wedge w) = u+(z\wedge w)$ that hold in any Riesz space.\\
	$R4)$ From Eq~\eqref{eqinf} we have
\[
i_\R(\omega)-\lambda = \min\{\pi_1(\Rj_\jg(\omega))-\lambda, \pi_2(\Rj_\jg(\omega))-\lambda\}.
\]
	Then $i_\R(\omega)-\lambda>0$ if and only if $\pi_1(\Rj_\jg(\omega))-\lambda>0 $ and $\pi_2(\Rj_\jg(\omega))-\lambda>0 $, i.e. if and only if $\Rj_\jg(\omega)\gg \lambda$ in $\D$.
\end{proof}

\subsection{Bicomplex Riemann zeta function}

The bicomplex Riemann zeta function is defined in~\cite{ZET} as the sum of the convergent series
\begin{equation*}
 \zeta(\omega): = \displaystyle\displaystyle \sum_{n = 1}^\infty \frac{1}{n^{\omega}}
\end{equation*}
on the open set
 $U : = \{\omega \in \B:\Rj(\mathcal{P}_1(\omega))>1 \mbox{~and~} \Rj (\mathcal{P}_2(\omega))>1 \}$. It's clear that $\omega \in U$ if and only if $\Rj_\jg(\omega)\gg 0$. Following our notations one has
 \begin{equation}\label{zetadomain}
   U = \{\omega \in \B: i_\R(\omega)>1\}.
 \end{equation}
Let us recall the main results of~\cite{ZET}.
\begin{itemize}
   \item The bicomplex Euler product formula :
   \begin{displaymath}
    \zeta(\omega) = \prod_{p \in \mathbb{P}} \dfrac{1}{1-\frac{1}{p^{\omega}}} \quad \mbox{for all\quad } \omega \in U.
   \end{displaymath}
   Where $\mathbb{P} $ is the set of all prime numbers.
   \item The analytic continuation of the bicomplex zeta Riemann function to the connected open set $1+ \B_\ast$ is defined by
   \begin{equation}\label{zetaprolong}
    \zeta(\omega) = \zeta(\mathcal{P}_1(\omega)) \eg_1+\zeta(\mathcal{P}_2(\omega)) \eg_2.
   \end{equation}
   \item The set of the trivial zeros for the the bicomplex Riemann zeta function is given by
   \begin{equation*}
    \mathcal{O} = \{\omega \in \B;\, \omega = (-n-p)+\jg (-n+p):\, n, p \in \N\setminus \{0\}\}.
   \end{equation*}
   Note that using the notion of hyperbolic integers introduced in definition~\ref{hypint} one can see that the trivial zeros of the bicomplex Riemann zeta function are
   \begin{equation*}
        \mathcal{O} = \{\nu = -2 h : h \in \ZG,\, h\geq 1 \}.
       \end{equation*}

   \item
   Riemann hypothesis (RH) is generalized to a bicomplex Riemann hypothesis ($\B$RH). It is shown that ($\B$RH) is equivalent to (RH).
 \end{itemize}

In the sequel we will use the Riesz subnorm $\|.\|_\jg$ defined by equation~\ref{SubN}.
\begin{theorem}\label{absconv}

\begin{itemize}
  \item [$(\rm{i})$] The series $ \displaystyle\sum_{n = 1}^{\infty} \dfrac {1}{n^{\omega}}$ is absolutely convergent if and only if $ \omega \in U$. Moreover,
    \begin{displaymath}
   \displaystyle \Big \|\sum_{n = 1}^{\infty} \dfrac {1}{n^{\omega}}\Big\|_{\jg}\leq \zeta[i_\R(\omega)].
  \end{displaymath}
  \item [$(\rm{ii})$] The series $ \displaystyle\sum_{n = 1}^{\infty}\dfrac{1}{n^{\omega}}$ is uniformly convergent on each compact subset of $U$.
\end{itemize}
 \end{theorem}
 \begin{proof}

   $(\rm{i})$ Let $ \omega \in \B$ and an integer $n\geq 1$. We have $\Ar_\D(n) = 0$ then Eq~\eqref{bicompawer} implies that $\dfrac{1}{n^{\omega}} = e^{-\omega \ln (n)}$. Therefore $ \Big| \dfrac{1}{n^{\omega}}\Big|_\jg = e^{-\Rj_\jg(\omega).\ln(n)}$. Since the hyperbolic exponential preserves lattice and conjugation operations (\ref{expop}) then by the identity $-z \vee -w = -(z\wedge w)$ we obtain
    \begin{equation*}
       \Big\| \dfrac{1}{n^{\omega}}\Big\|_\jg = e^{s_\R(- \ln (n)\omega)}.
    \end{equation*}
 But $s_\R(- \ln (n)\omega) = - \ln (n) i_\R(\omega)$ (by $R3)$). Hence
   \begin{equation*}
   \Big\| \dfrac{1}{n^{\omega}}\Big\|_\jg = \dfrac{1}{n^{i_\R(\omega)}}.
   \end{equation*}
This proves that the series $ \displaystyle\sum_{n = 1}^{\infty}\dfrac{1}{n^{\omega}}$ is absolutely convergent if and only if $i_\R(\omega)> 1 $, i.e. (by Eq~\eqref{zetadomain}) if and only if $ \omega \in U$. Moreover,
        \begin{displaymath}
      \Big \|\sum_{n = 1}^{\infty}\dfrac{1}{n^{\omega}} \Big\|_\jg \leq \sum_{n = 1}^{\infty} \Big \|\dfrac{1}{n^{\omega}} \Big \|_\jg = \zeta[i_\R(\omega)].
    \end{displaymath}
       $(\rm{ii})$ Let $K$ be a compact subset of $U$. Then, $i_\R(K)$ is a compact of
   $i_\R(U) = i_\R [i_\R^{-1}[(1,\infty)]] = (1,\infty)$, since $i_\R: \B\longrightarrow \R$ is a surjective and continuous function. Which implies that $ i_\R(K)\subset [1+\alpha,\infty)$ for some real $\alpha>0$ and then
\[
K\subset i_\R^{-1} i_\R (K)\subset i_\R^{-1}[ \{[1+\alpha,\infty)\}] = U_\alpha.
\]
 Therefore
 \begin{align*}
 a_n = \sup_{\omega \in K} \Big\| \frac{1}{n^{\omega}}\Big\|_\jg 
 &\leq \sup_{\omega \in U_\alpha} \Big\| \frac{1}{n^{\omega}}\Big\|_\jg\\
 &\leq \sup_{\omega \in U_\alpha} \frac{1}{n^{i_\R(\omega)}}\\
 &\leq \frac{1}{n^{1+\alpha}}.
 \end{align*}
 This means that $\displaystyle\sum_{n = 1}^{\infty} a_n $ is convergent and hence $ \displaystyle\sum_{n = 1}^{\infty}\dfrac{1}{n^{\omega}}$ is uniformly convergent in~$K$.
\end{proof}
\subsection{Bicomplex gamma function}

The bicomplex gamma function was first introduced in~\cite{GOYAL} by an Euler product formula using the idempotent decomposition on $\B $. Therefore, it is defined as a pair of complex gamma functions. Our approach is novel and extend classical proof to bicomplex numbers as a Banach algebra equipped with a Riesz subnorm.
Contrary to~\cite{GOYAL} we introduce the bicomplex gamma function as an absolute convergent integral.
 \begin{proposition}[Bicomplex gamma function]
 	Let $\omega \in \B$ be such that $\Rj_\jg(\omega)\gg 0$, then the integral
 	
\begin{equation*}
 	\Gamma(\omega) = \int_0^\infty e^{-t} t^{\omega-1} dt
\end{equation*}
 	 called bicomplex gamma function, is absolute convergent.
 \end{proposition}
\begin{proof}
	We have that $|e^{-t} t^{\omega-1}|_\jg = e^{-t} e^{\ln(t) (\Rj_\jg(\omega)-1)}$ and then
\[
	\| e^{-t} t^{\omega-1}\|_\jg = e^{-t} e^{s_\R\left[ \ln(t) (\omega-1) \right]}.
\]
	From $R3)$ we have $s_\R\left[ \ln(t) (\omega-1) \right] = {(\ln t)}^+ (s_\R(\omega)-1)-{(\ln t)}^{-} (i_\R(\omega)-1) $. Thus
	
\begin{equation*}
	\| e^{-t} t^{\omega-1}\|_\jg\sim t^{i_\R(\omega)-1} \quad \mbox{as~} t\longrightarrow 0^+
\end{equation*}
	and
	\begin{equation*}
	\| e^{-t} t^{\omega-1}\|_\jg\sim e^{-t} t^{s_\R(\omega)-1} \quad \mbox{as~} t\longrightarrow \infty
\end{equation*}
	Assume that $\Rj_\jg(\omega)\gg0$. So, $i_\R(\omega),s_\R(\omega) >0$, since $\D^+_\ast$ is closed under conjugation and lattice operations. It follows that the integral defining $\Gamma(\omega)$ is absolute convergent.
\end{proof}
The next result extend Mellin integral to bicomplex numbers
\begin{theorem}[Bicomplex Mellin Integral ]
	For $\Rj_\jg(\omega)\gg 1$ we have the following integral representation of bicomplex zeta function
\begin{equation*}
	\zeta(\omega) \Gamma (\omega) = \int_0^\infty
 \frac{t^{\omega-1}}{e^t -1} dt.
\end{equation*}
\end{theorem}
\begin{proof}
	Put $t = n s$ in the integral representation of $\Gamma$, one gets
	
\begin{equation*}
	\zeta(\omega)\Gamma (\omega) = \displaystyle \lim_{N\rightarrow\infty} \Phi_N(\omega),
\end{equation*}
	where $\Phi_N(\omega) = \displaystyle \sum_{n = 1}^N \int_0^\infty s^{\omega-1} e^{-n s} ds$. We have
	
\begin{equation*}
	\Phi_N(\omega) = \int_0^\infty \frac{s^{\omega-1}}{e^s -1} ds-\int_0^\infty \frac{s^{\omega-1}}{e^s -1} e^{-Ns} ds.
\end{equation*}
	The above integrals are absolutely convergent for $\Rj_\jg(\omega)\gg 1$. Since $i_\R(\omega) > 1$, $s_\R(\omega)> 1$, and for every integer $N\geq 0$ we have
\[
\Big\|\frac{s^{\omega-1}}{e^s -1} e^{-Ns}\Big\|_\jg	= \left\{
\begin{array}
{rl} \frac{s^{i_\R(\omega)-1}}{e^s -1} e^{-Ns} & \mbox{if~} 0< s\leq 1\\
\\
\frac{s^{s_\R(\omega)-1}}{e^s -1} e^{-Ns}		& \mbox{if~} s\geq 1.
\end{array}
\right.
\]
Set
\begin{equation*}
	I_N(\omega) = \int_0^\infty \frac{s^{\omega-1}}{e^s -1} e^{-Ns} ds = \int_0^\infty s^{\omega-2}\frac{s} {e^s -1} e^{-Ns} ds.
\end{equation*}
We have
\begin{align*}
\|I_N(\omega)\|_\jg\leq \int_0^\infty \| s^{\omega-2} e^{-N s}\|_\jg ds = \frac{1}{N}\int_0^\infty e^{-t} \|{(\frac{t}{N})}^{\omega-2}\|_\jg dt.
\end{align*}
By $\|{(\frac{t}{N})}^{\omega-2}\|_\jg = {(\frac{t}{N})}^{i_\R(\omega)-2}$ for $0<t\leq 1$ and $\|{(\frac{t}{N})}^{\omega-2}\|_\jg = {(\frac{t}{N})}^{s_\R(\omega)-2}$ for $ t\geq1$ we obtain
\begin{align*}
\frac{1}{N}\int_0^\infty e^{-t} \|{(\frac{t}{N})}^{\omega-2}\|_\jg dt 
&= \frac{1}{N} \left(\int_0^1 e^{-t} {(\frac{t}{N})}^{i_\R(\omega)-2} dt+ \int_1^\infty e^{-t} {(\frac{t}{N})}^{s_\R(\omega)-2} dt\right)\\
&\leq \frac{1}{N^{i_\R(\omega)-1}} \Gamma(i_\R(\omega)-1)+\frac{1}{N^{s_\R(\omega)-1}} \Gamma(s_\R(\omega)-1).
\end{align*}

	Therefore, $\displaystyle \lim_{N\rightarrow\infty} I_N(\omega) = 0$ and hence
\begin{equation*}
\zeta(\omega)\Gamma (\omega) = \displaystyle \lim_{N\rightarrow\infty} \Phi_N(\omega) = \int_0^\infty
\frac{t^{\omega-1}}{e^t -1} dt.
\qedhere
\end{equation*}
\end{proof}

\begin{theorem}[Analytic continuation]
By analytic continuation, the bicomplex gam\-ma function is $\B$-holomorphic on
\[
\Omega_-: = \{\omega \in \B;\,\mathcal{P}_1(\omega) \notin \Z^- \mbox{~and~}\mathcal{P}_2(\omega) \notin \Z^- \}.
\]
\end{theorem}
\begin{proof}
Using idempotent decomposition of $\B$, for each $\omega = z_1 \eg_1+z_2 \eg_2 $ we have $e^{-t} t^{\omega-1} = e^{-t} t^{z_1} \eg_1+e^{-t} t^{z_2-1} \eg_2$. Thus, for $\Rj_\jg(\omega)\gg 0$ we have $\Rj(z_1)>0$ and $\Rj(z_2)>0$. Therefore
\begin{equation}\label{decompgamma}
\Gamma(z_1\eg_2+z_2 \eg_2) = \Gamma(z_1)\eg_1+\Gamma(z_2)\eg_2.
\end{equation}
Knowing that the classical complex gamma function is extended, by analytic continuation, to a holomorphic function on
$\C \setminus \Z^-$,
the representation~\eqref{decompgamma} allows us to extend the bicomplex gamma function $\Gamma$ to a $\B$-holomorphic function on $\Omega$ as follows
\begin{equation}\label{prologamma}
\Gamma(\omega) = \Gamma(\mathcal{P}_1(\omega)) \eg_1+\Gamma(\mathcal{P}_2(\omega)) \eg_2.
\qedhere
\end{equation}
\end{proof}
\begin{corollary}[Weierstrass formula]
For $\omega \in \Omega$ one has
\begin{equation*}
\displaystyle \Gamma(\omega) = \frac{e^{-\gamma \omega}}{\omega \displaystyle \prod_{n = 1}^{\infty} \left(1+ \frac{\omega}{n}\right) e^{\frac{-\omega}{n}}}.
\end{equation*}
\end{corollary}
\begin{proof}
From~\eqref{prologamma} we have that $\mathcal{P}_k(\Gamma(\omega)) = \Gamma(\mathcal{P}_k(\omega))$, $(k = 1,2)$. Thus, applying the classical Weierstrass formula to
$\mathcal{P}_k(\omega) \notin \Z^-,\, (k = 1,2), $ one gets the above representation of $\Gamma(\omega).$
\end{proof}

Now we can achieve the proof of Theorem~\ref{main}.
According to circular functions defined on $\D$ (see section 4.1) one
can define for a
given bicomplex number $\omega$, $\sin \omega$ as
\begin{equation}\label{sine}
\sin \omega = \displaystyle \sum_{n = 0}^\infty {(-1)}^n \frac{\omega^{2n+1}}{(2n+1)!} = \sin (\mathcal{P}_1(\omega)) \eg_1+\sin (\mathcal{P}_2(\omega))\eg_2.
\end{equation}
Note that $\sin (\pi\omega)$ is invertible if and only if
\[
\omega \in \Omega: = \{\omega \in \B;\,\mathcal{P}_1(\omega) \notin \Z \mbox{~and~}\mathcal{P}_2(\omega) \notin \Z\}.
\]
Then, Theorem~\ref{main} follows immediately from usual complex functional equations, using the idempotent representations of the bicomplex zeta, gamma and sine functions given by~\eqref{zetaprolong},~\eqref{prologamma} and~\eqref{sine}.

\EditInfo{%
    June 18, 2020}{%
    August 09, 2020}{%
    Valentin Ovsienko
    }

\end{paper}